\documentclass[10pt]{article}
\usepackage[utf8]{inputenc}

\usepackage[letterpaper, margin=1in]{geometry}
\usepackage{parskip}

\usepackage{sectsty}
\sectionfont{\fontsize{12}{14.4}\selectfont}

\usepackage{titlesec}
\titlespacing*{\section}{0pt}{10pt}{0pt}

\usepackage{caption}
\usepackage{setspace}
\usepackage[T1]{fontenc}
\usepackage[english]{babel}

\usepackage{ifpdf,newtxtext,newtxmath} 

\usepackage{array,graphicx,dcolumn,multirow,hevea,abstract,hanging,fancyhdr}
\usepackage{algorithmic}
\usepackage[ruled,vlined]{algorithm2e}
\usepackage {tikz}
\usetikzlibrary {positioning}
\def\T{\mathcal{T}}
\def\E{\mathcal{E}}
\def\N{\mathcal{N}}
\def\Gv{\mathcal{G}}
\def\D{\mathcal{D}}

\def\M{\mathcal{M}}

\def\K{\mathcal{K}}

\newcommand{\jheadfirstpage}{\itshape Proceedings of the 2021 IISE Annual Conference \\
A. Ghate, K. Krishnaiyer, K. Paynabar, eds.}
\newcommand{\jheadrestpages}{}

\fancypagestyle{firstpage}{
 \lhead{\protect\fontsize{10}{12}\jheadfirstpage}
 \rhead{}
}

\usepackage[labelfont=sc,textfont=sf]{caption}
\usepackage[hyperfootnotes=false,breaklinks=true]{hyperref} 
\usepackage[hyphenbreaks]{breakurl}

\usepackage{booktabs} 
\newcolumntype{L}[1]{>{\raggedright\arraybackslash }p{#1}} 
\newcolumntype{C}[1]{>{\centering\arraybackslash }p{#1}}
\newcolumntype{R}[1]{>{\raggedleft\arraybackslash }p{#1}}
\newcolumntype{d}[1]{D{.}{.}{#1}} 

\let\tempone\itemize
\let\temptwo\enditemize
\let\tempthree\enumerate
\let\tempfour\endenumerate
\renewenvironment{itemize}{\tempone\setlength{\itemsep}{0pt}}{\temptwo}
\renewenvironment{enumerate}{\tempthree\setlength{\itemsep}{0pt}}{\tempfour}

\newcommand{\zred}{\textcolor{black}}

\title{\fontsize{16}{19.2} \textbf{Distributionally Robust Optimal Power Flow with Uncertain Renewable Energy Output}}

\author{
\fontsize{12}{14.4}
\textbf{Jia Yang} \\ 
\fontsize{12}{14.4}
\textbf{Power Costs, Inc.} \\
\fontsize{12}{14.4}
\textbf{Norman, OK 73072, USA} \\
\\
\fontsize{12}{14.4}
\textbf{Jun Song and Chaoyue Zhao} \\
\fontsize{12}{14.4}
\textbf{Department of Industrial and Systems Engineering, University of Washington} \\
\fontsize{12}{14.4}
\textbf{Seattle, WA 98195, USA}
}

\date{} 

\begin{document} 
\singlespacing

\maketitle
\thispagestyle{firstpage}

\section*{\centering Abstract}
\vskip 10pt
Optimal power flow (OPF) is an important tool for Independent System Operators (ISOs) to deal with the power generation management. With the increasing penetration of renewable energy into power grids, challenges arise in tackling the OPF problem due to the intermittent nature of renewable energy output. To address these \zred{challenges}, we develop a multi-stage distributionally robust
approach for the direct-current optimal power flow (DC-OPF) problem to
minimize total generation cost under renewable energy uncertainty. In our model, we assume the renewable energy output follows an ambiguous distribution that can be characterized by a confidence set. By utilizing the revealed data sequentially, the proposed approach can provide a reliable and robust optimal OPF decision without restricting the renewable energy output distribution to any particular distribution class. The computational results also verify the effectiveness
of our approach to reduce the \zred{conservativeness} and meanwhile maintain the reliability. 

\section*{Keywords} 
Optimal Power Flow, Renewable Energy Output Uncertainty, Distributionally Robust Dynamic Programming

\section{Introduction}
Optimal power flow (OPF), as an essential component of power system operations and management, is aimed at controlling the power flow without violating any operational restriction. As the penetration of renewable energy output has been increasing dramatically, how to efficiently and securely address the OPF problem remains challenging due to the intermittency and uncertainty that the renewable energy brings to the power system operations. \\

The most common approach to formulate the OPF problem with renewable energy uncertainties is stochastic OPF (see, e.g., \cite{perninge2012stochastic,liang2013two,xia2016probabilistic}, among others), in which the uncertain renewable energy output is characterized by a finite number of scenarios or assumed to follow a particular probability distribution. Another traditional approach to address the OPF problem with renewable energy uncertainties is the robust optimization approach (see, e.g., \cite{c13,cc40,cc41}, among others), which constructs a solution that is optimal for the worst-case realization of renewable energy output in a predefined uncertainty set. However, both the stochastic OPF models and the robust OPF models have disadvantages in practice. For the stochastic approach, the obtained solution can be sub-optimal or biased if the assumption of distribution is inaccurate. For the robust approach, since the optimal solution is based on the worst-case scenario, which happens rarely, the optimal objective and solution are usually over conservative and cost ineffective. In addition, \zred{robust optimization} cannot fully utilize data information, since it only requires minimal information (such as lower/upper bounds of the uncertain parameter) to generate the uncertainty set. \\

In order to bridge the gap between the stochastic optimization approach and the robust optimization approach, and \zred{to} utilize data information as much as possible, in this paper, we explore distributionally robust optimization (DRO) approaches to solve the OPF problem. Instead of considering a pre-assumed distribution like stochastic OPF or the worst case scenario like robust OPF, we assume that the distribution is unknown but can be characterized by a confidence set constructed by learning from the historical data. The objective of distributionally robust OPF problem is to minimize the total cost under the worst-case probability distribution of renewable energy output within \zred{a predefined} confidence set. DRO-based
approaches have recently received attention from power system
specialists for various problems such as unit commitment \cite{zhao2017distributionally}, transmission expansion planning and hardening \cite{bagheri2016data},
and reserve scheduling \cite{wei2015distributionally}. To the best of our knowledge, DRO is only used in \cite{guo2018data} to address the multi-stage OPF problem. However, \cite{guo2018data} adopts a strong assumption: it assumes that the power balance constraints are linear functions of the ambiguous distribution, which is not practical for real OPF problems. In this paper, we do not pre-assume any relationship between the OPF control policy and the random renewable energy output. \\

The rest of this paper is organized as follows. In Section 2, we develop a distributionally robust dynamic direct-current optimal power flow model under the renewable energy output uncertainty. Then, in Section 3, we reformulate the problem as a tractable second order cone problem (SOCP), and propose the corresponding algorithm to solve the reformulated SOCP. In Section 4, we conduct experiments on a 20-bus system to show the effectiveness of the \zred{proposed} model. In Section 5, we conclude our work.

\section{Model Formulation}
We first describe the detailed formulation of the multi-stage distributionally robust direct-current optimal power flow problem. We use $\T$, $\N$, $\Gv$, $\E$, $\M$, $\K$ to represent sets of time periods, buses, generators, transmission lines, renewable energy output scenarios, and interpolation points of generation amount respectively. The parameters of each generator $i$ include minimal power output $P_{min}^i$, maximal power output $P_{max}^i$, ramp-up limit $R^{u}_i$ and ramp-down limit $R^{d}_i$. For each transmission line $(m,n) \in \E$, $X_{mn}$ represents the reactance and $L_{mn}$ represents the transmission capacity. For each bus $n$, $\theta_{min}^n$ and $\theta_{max}^n$ represent the minimal and maximal values of the phase angle respectively. For generator $i$ at bus $n$ in period $t$, $D_{nt}$ represents the demand, $\xi_{nt}$ represents the renewable generation, $C_{t}^i(\cdot)$ represents the fuel cost function, and $P_n^+$ and $P_n^-$ represent the unit penalty cost for load shedding and renewable energy curtailment respectively. The decision variables include generation level ($x_{t}^i$) of generator $i$ in period $t$, power flow ($p_{mn}^t$) on transmission line ($m,n$) in period $t$, load shedding ($q_{nt}^+$) at bus $n$ in period $t$, renewable energy curtailment ($q_{nt}^-$) at bus $n$ in period $t$, phase angle ($\theta_{nt}$) at bus $n$ in period $t$, and generator $i$'s ramping amount ($\Delta x_{t}^i$) at period $t$. In addition, $\theta_t(ref)$ represents the phase angle of the reference bus. The detailed formulation is shown as follows: 
\begin{eqnarray}
\min \quad
\sum_{t\in\T} \sum_{i\in\Gv} C_{t}^i(x_{t}^i)+ \sum_{t\in\T} \sum_{n\in\N} (P_n^+q_{nt}^+ +P_n^-q_{nt}^-)\label{obj}\\
s.t.\hspace{6mm}\underset{i\in\Gv_n}{\sum x_{t}^i}+\underset{m\in\E(.,n)}{\sum}{p_{mn}^t}-\underset{m\in\E(n,.)}{\sum}{p_{mn}^t} + 
q_{nt}^+- q_{nt}^-& \nonumber \\ 
= D_{nt}-\xi_{nt},\hspace{5mm}
\forall t\in\T,\  \forall n\in\N,\ \forall i\in\Gv_n, \label{balance}\\
P_i^{min} \leq x_{t}^i \leq P_i^{max}, \ \  \forall i\in\Gv,\ \forall t\in\T, \label{capacity}\\
\Delta x_{t}^i = x_t^i - x_{t-1}^i, \forall t\in \T \setminus \{1\}, \ \  \Delta x_{1}^i = x_1^i, \ \ \forall i\in\Gv, \label{rampup}\\
\Delta x_{t}^i \leq R_i^u, \ \ 
-\Delta x_{t}^i \leq R_i^d,\ \ \forall i\in\Gv,\  \forall t\in\T, \label{rampdown}\\
\theta_t(ref) = 0, \ \  \forall t\in\T, \label{ref}\\
\theta_{min}^n\leq \theta_{nt}\leq \theta_{max}^n,
\ \  \forall t\in \T,\ \forall n\in \N, \label{angle}\\
-L_{mn}\leq p_{mn}^t\leq L_{mn}, \ \  \forall (m,n)\in\E,\  \forall t \in \T,\label{transcap}\\
(\theta_{mt}-\theta_{nt})/X_{mn} = p_{mn}^t,\ \  \forall (m,n)\in \E,\  \forall t\in\T. \label{angletrans}
\end{eqnarray} 

In the above formulation, the objective function is to minimize the total cost \zred{which} includes the fuel costs and the penalty costs; constraint \eqref{balance} represents the energy balance constraint; constraint \eqref{capacity} \zred{specifies} generation capacity limits for each generator. The ramp-up and ramp-down restrictions are described in \eqref{rampup} and \eqref{rampdown}, \zred{and} constraint \eqref{ref} fixes the phase angle of the reference bus to zero. The restrictions of phase angle and transmission capacity are described in constraints \eqref{angle} and \eqref{transcap}, and constraint \eqref{angletrans} describes the relationship between the power flow and the phase angle for each transmission line. Note here we just assume that the renewable energy is uncertain, but this proposed framework can be used to \zred{incorporate} both demand and renewable uncertainties. In this case, we can define the random parameter $\xi$ as the net load, which is the load demand minus the renewable energy output. \\

Since the OPF problem is inherently sequential, which means the information of the uncertain renewable output $\xi$ is constantly revealed as the time progresses, and dispatch decisions at each period are made after knowing the realization of uncertain parameters up to that period, we use \zred{a} stochastic dynamic programming framework to solve the DC-OPF problem. In addition, we assume that the distribution of renewable output is unknown but \zred{lies} within \zred{the} confidence set $\Psi$ defined as follows:
\begin{align}
&\Psi := \{p:\ \sum_{j \in \M}p_j = 1, \ \  \sum_{j\in \M}(p_j-q_j)^2/p_j\leq \gamma\}.\label{uncertain}
\end{align} 

Here we assume that there are a total number of $|\M|$ scenarios of the renewable energy output. And for each scenario $j$, we use $p_j$ and $q_j$ to represent the true probability and the reference probability learned from data respectively. The confidence set $\Psi$ is constructed by using $\chi^2$ divergence but the framework also works \zred{with} a general distribution metric. $\gamma$ is the tolerance level of the distance between the reference distribution and the true distribution, and it can be used to control the robustness of this framework: the larger $\gamma$ value is, the higher level of robustness the system has. \\

In this paper, we propose a distributionally robust dynamic OPF (DRD-OPF) framework by considering the worst-case distribution within the confidence set $\Psi$, and minimizing the total cost under the worst-case distribution. The Bellman equation can be expressed as follows:
\begin{align}
f_t(x_t,\xi_t) = &\min_{p_t,\theta_t,\Delta x_t \in \D_t}\sum_{i\in\Gv}C_{t}^i(x_{t}^i)+  \sum_{n\in\N} (P_n^+q_{nt}^+ +P_n^-q_{nt}^-)+\max_{p\in\Psi}\sum_{j\in\M}p_jf_{t+1}(x_{t+1}^{j},\xi_{t+1}^{j})\label{MSobj}\\
s.t.\hspace{6mm}&
\D_t:=\{\text{Constraints}\ \eqref{balance}-\eqref{angletrans} \text{ for time period } t\}.
\end{align} 

In the above formulation, the objective function in \eqref{MSobj} is to minimize the total cost from time period $t$ to $T$ with terminate condition $f_{T+1}(.) = 0$. The generation amount $x$ is considered as the decision-dependent (endogenous) state variable, and the renewable energy output $\xi$ is considered as the decision-independent (exogenous) state variable. The ramping amount, the power flow and the phase angle are considered as action variables.

\section{Solution Approach}
To solve the above formulation, we employ the method proposed by Hanasusanto and Kuhn \cite{cc28}. We first reformulate the DRD-OPF into a second order conic program by dualizing the second-stage maximization problem in \eqref{MSobj}: 
\begin{align}
f_t(x_{t},\xi_t)=&\min_{p_t,\theta_t, \Delta x_t \in \D_t}\ \sum_{i\in\Gv}
C_{t}^i(x_{t}^i)+ \sum_{n\in\N} (P_n^+q_{nt}^+ +P_n^-q_{nt}^-)+
\beta\gamma-\sigma-2q_t^T\textit{y}+2\beta\\
\mbox{(Dual)} \ \ s.t. \ \ &\sqrt{4y_j^2+(z_j+\sigma)^2}\leq2\beta-z_j-\sigma,\hspace{3mm}\forall j\in \M,\\
&f_{t+1}(x_{t+1}^j, \xi_{t+1}^j)\leq z_j, \hspace{3mm}  \forall j \in \M,\label{costf}\\
&z_j + \sigma  \leq\beta,\ \ \forall j \in \M,\\
&\sigma \in\mathbb{R}, \ \ \beta  \in \mathbb{R}_+, \ \
\textbf{z},\textbf{y}  \in \mathbb{R}^{|\M|},
\end{align}
where $\beta$ and $\sigma$ are dual variables for constraints in the confidence set $\Psi$. \\

Note that in the above formulation, the left hand side in constraint \eqref{costf} represents the total optimal cost from period $t+1$ to $T$. To approximate the cost function $f$, we follow the idea in \cite{cc28} to use a quadratic function  $x^TAx+2B^Tx+C$, which minimizes the total square error between $f$ and the quadratic function at a series of interpolation points of generation amount $x_t^{k}$, for $k \in \K$ and $t \in \T$.  The formulation is shown as follows:
\begin{align}
&\min \quad \sum_{k\in \K} ((x_{t}^{k})^TA_jx_{t}^k+2B_{j}^Tx_{t}^{k}+C_j-
\hat{f}_{t}^{j,k})^2\\
(\mbox{QAP})\ \  &s.t.\ \ \ A_j \in \mathbb{S}^{|\Gv|},\ A_j \succeq 0,\  B_j \in \mathbb{R}^{|\Gv|},\ C_j \in\mathbb{R}, \ \forall j \in \M.
\end{align} 

We solve the dual problem (Dual) backwardly from time period $t = T$. When $t = T$, we can obtain the optimal solution $\hat{f}_T^{j,k}$ based on the input $x_{T}^{k}$ and $\xi_T^j$, for $j \in \M, k \in \K$. In this period, constraint \eqref{costf} can be simplified as $z_j\geq0$ for $j \in \M$, 
based on the initialization of the algorithm. Then, we can obtain the optimal value of (Dual) for every  generation amount $x_t^k$ for each scenario $j$ at time period $t = T$.  Then we solve the semidefinite program (QAP) to obtain the parameters $A_j, B_j$ and $C_j$. Once the semidefinite program is computed for all $j \in \M$, we update the constraint \eqref{costf} using the quadratic function $x^TA_jx+2B_j^Tx+C_j$. The algorithm proceeds iteratively until period \zred{$t = 1$}. 


\section{Case Study}
In this section, we perform a case study on a 20-bus system to show the system performance of the proposed DRD-OPF model. The 20-bus system contains 6 thermal generators, 20 loads and 45 transmission lines. All the experiments are implemented in C++ with CPLEX 12.7 and Python with MOSEK. The time interval for all experiments is 1 hour and the time horizon is 24 hours. \\

We set the penalty cost for load-shedding and renewable energy curtailment as \$1,500/MW. In addition, we partition the range of generation amount for each generation into 3 pieces and set the central point of each piece as the interpolation point. Therefore we consider $3^6 = 729$ interpolation points in total. In order to show the system performance under different system robustness levels, we test five cases of robustness parameter $\gamma$: 0, 1, 3, 5, and 10, and report the total costs in \textbf{Table \ref{mylabel}}. From \textbf{Table \ref{mylabel}}, we can observe that the total cost increases as the robust parameter $\gamma$ increases, as the algorithm becomes more conservative and cost effectiveness is sacrificed to maintain a higher level of robustness. Note here when $\gamma$ equals to 0, the model becomes a traditional multi-stage stochastic program (MSSP), and the proposed model is more conservative compared with the MSSP as it yields more costs. \\


\begin{table}[h!]
\caption{\textrm{Total Cost vs. Robustness Level}}
\centering
\begin{tabular}{|c|c|c|c|c|c|}
\hline
\textit{$\gamma$} & 0 & 1 & 3& 5 & 10\\
\hline
Total Cost (\$) &253995  &262028 &263883 &264781 &269748
\\ 
\hline
\end{tabular}
\label{mylabel}
\end{table}\par

We further compare the performance of our proposed model with the traditional MSSP (i.e., $\gamma = 0$) via simulation. First, we generate $500$ samples of renewable energy output and obtain the reference distribution as the histogram of the $500$ data samples. Then, we obtain optimal generation scheduling solutions for both the proposed approach and the MSSP with the reference distribution. Then we fix the optimal generation scheduling for each approach, and simulate another 1,000 samples to obtain the simulated total costs under the previously fixed generation scheduling for each approach. We report the mean value, the standard deviation and the 90th percentile value of the results for each approach in \textbf{Table \ref{table2}}. \\

\begin{table}[h]
\caption{\textrm{Comparison with Traditional MSSP}}
\centering
\label{my-label}
\begin{tabular}{|c|c|c|c|}
\hline
$\gamma$  & Mean  &  Std. Dev.& 90th  percentile  \\
\hline

0   & 259864.5         &   1261.9 &261478.3               \\
1   & 260748.1     &    1208.7 & 262197.6               \\
3   & 261041.3         &  970.64 &261504.5                  \\
5   & 259072     &  744.2 &260403.2                  \\
10   & 257433.3       &  544.7  & 259982.6          \\
\hline
\end{tabular}
\label{table2}
\end{table}\par

From \textbf{Table \ref{table2}}, we can observe that, if the robustness level of the system is high, i.e., $\gamma = 5$ and $10$, the proposed approach outperforms the MSSP approach since it provides a \zred{more} robust and reliable solution. However, when the robustness level is low, which means we have relatively accurate information of renewable energy output, the MSSP approach is better since it is a risk-neutral approach and is not too conservative as compared with the proposed DRD-OPF approach. 

\section{Conclusion}
In this paper, we develop a distributionally robust dynamic programming approach for the direct-current optimal power flow problem under renewable energy uncertainty. Our proposed approach does not require any assumptions of the probability distribution of renewable energy output, as the traditional stochastic optimization approach does.  Instead, we construct a confidence set of distribution based on historical data information, and consider the worst case distribution within the set. The advantage of our model is that we can guarantee the robustness of the model without ignoring the stochastic nature of uncertain renewable energy, and \zred{we can} utilize the revealed data sequentially and efficiently. Meanwhile, the conservativeness of our proposed approach is adjustable based on system operator's preference. Furthermore, this proposed distributionally robust dynamic programming framework not only helps solve the OPF problem, but
also helps solve other power system problems such as generation
investment, transmission planning and contingency analysis. 

\section*{Acknowledgement}
This research was supported by National Science Foundation Grant \#1662589. This manuscript has benefited from the comments and suggestions of anonymous reviewers.

\bibliographystyle{IEEEtran}
\bibliography{template}

\end{document}